\documentclass[11pt]{amsart}
\usepackage{amsmath,latexsym,amssymb}
\usepackage{mathrsfs}
\usepackage[english]{babel}
\usepackage[applemac]{inputenc}
\usepackage{amsthm}
\usepackage[all]{xy}

\renewcommand{\>}{\rangle}
\newcommand{\<}{\langle}
\newtheorem{definition}{Definition}[section]
\newtheorem{theorem}{Theorem}[section]
\newtheorem{proposition}[theorem]{Proposition}
\newtheorem{lemma}[theorem]{Lemma}
\newtheorem{remark}{Remark}[section]
\newtheorem{corollary}[theorem]{Corollary}
\newtheorem{question}{Question}
\newtheorem{example}{Example}
\newtheorem{notation}{Notation}[section]
\def\bem{\begin{remark}\upshape}
\def\ebem{\end{remark}}
\def\nota{\begin{notation}\upshape}
\def\enota{\end{notation}}
\def\defn{\begin{definition}\upshape}
\def\edefn{\end{definition}}
\def\thm{\begin{theorem}}
\def\ethm{\end{theorem}}
\def\lmm{\begin{lemma}}
\def\elmm{\end{lemma}}
\def\qed{\hfill$\quad\Box$}
\def\pr{\par\noindent{\em Proof: }}

\def\kor{\begin{corollary}}
\def\ekor{\end{corollary}}
\def\frag{\begin{question}\upshape}
\def\efrag{\end{question}}
\def\prop{\begin{proposition}}
\def\eprop{\end{proposition}}

\def\bsp{\begin{example}}
\def\ebsp{\end{example}}
\def\U{{\mathscr U}}

\def\Z{{\mathbb Z}}
\def\R{{ \mathbb R}}
\def\N{{\mathbb N}}
\def\L{{\mathcal L}}
\def\Fa{{\mathcal F}}
\def\P{{\mathcal P}}
\def\B{{\mathcal B}}
\def\A{{\mathcal A}}
\def\C{{\mathcal C}}

\def\D{{\mathcal D}}
\def\si{{\sigma}}
\title{Alternatives for pseudofinite groups}
\author{A. Ould Houcine}
\address{\hskip-\parindent
A. Ould Houcine\\
D\'epartement de math{\'e}matique, Le Pentagone\\
Universit{\'e} de Mons\\
20, place du Parc, B-7000 Mons, Belgium. Universit\'e de Lyon; Universit\'e Lyon 1; INSA de Lyon, F-69621; Ecole Centrale de Lyon; 
CNRS, UMR5208, Institut Camille Jordan, 43 blvd du 11 novembre 1918, F-69622 Villeurbanne-Cedex, France. }
 \email
{ould@math.univ-lyon1.fr}
\author{F.  Point}
\address{\hskip-\parindent
Fran\c coise Point\\
D\'epartement de math{\'e}matique, Le Pentagone\\
Universit{\'e} de Mons\\
20, place du Parc, B-7000 Mons, Belgium.}
\email {point@math.univ-diderot.fr}

\begin{document}

\begin{abstract} The famous Tits' alternative  states that a linear group either contains a nonabelian free  group or is soluble-by-(locally finite). We study in this paper similar alternatives in  pseudofinite groups. We show for instance that an $\aleph_{0}$-saturated pseudofinite group either contains a subsemigroup of rank $2$ or is nilpotent-by-(uniformly locally finite). We call a class of finite groups $G$ weakly of bounded rank if the radical $rad(G)$ has a bounded Pr\"ufer rank and the index of the sockel of $G/rad(G)$ is bounded. We show that an  $\aleph_{0}$-saturated pseudo-(finite weakly of bounded rank) group either  contains a nonabelian free group or is nilpotent-by-abelian-by-(uniformly locally finite). We also obtain some relations between  this kind of alternatives and amenability.
\end{abstract}
\maketitle



\section{Introduction}
A group $G$ (respectively a field $K$) is \textit{pseudofinite} if it is elementary equivalent to an  ultraproduct of finite groups (respectively of finite fields), equivalently if $G$ is a model of the theory of the class of finite groups (respectively of finite fields); that is any sentence true in $G$ is also true in some finite group. 
Note that one usually requires in addition the structure to be infinite, but it is convenient for us to allow a pseudofinite structure to be infinite. 

Infinite pseudofinite fields have been characterized algebraically by J. Ax \cite{Ax-pseudo} and he showed that the theory of all pseudofinite infinite fields is decidable.  Natural examples of pseudofinite groups are general linear groups over pseudofinite fields.  Pseudofinite simple groups have been investigated first by U. Felgner \cite{F}, then by J. Wilson \cite{W0} who showed that any pseudofinite simple group is elementarily equivalent to a Chevalley group (of twisted or untwisted type) over a pseudofinite field and it was later observed that it is even isomorphic to such a group \cite{Ry}. Pseudofinite groups with a theory satisfying various model-theoretic assumptions like stability, supersimplicity or the non independence property (NIP) have been studied \cite{MT, EJMR}; in another direction G. Sabbagh and A. Kh\'elif investigated finitely generated pseudofinite groups.

The Tits alternative \cite{tits} states that a linear group, i.e. a subgroup of the general linear group $GL(n,K)$ for some field $K$,  either contains a free nonabelian group or is soluble-by-(locally finite).   It is known that the Tits alternative holds for other classes of groups. For instance a subgroup of a hyperbolic group satisfies a strong form of the Tits alternative, namely it is either virtually cyclic or contains a nonabelian free group.  N. Ivanov \cite{Ivan}  and J.  McCarthy \cite{MacC}  have shown that mapping  
class groups of compact surfaces  satisfy the Tits  alternative and M. Bestvina, M. Feighn and M. Handel  \cite{F-Bes} showed that the alternative  holds for $Out(F_n)$  where $F_n$ is the free group of rank $n$. Note that when Tits alternative holds in a class of groups, then the following dichotomies hold for their finitely generated members: they have either polynomial or exponential growth; they are either amenable or contain a free nonabelian group. 
However, it is  well-known that groups which are non-amenable and without nonabelian free subgroups exist \cite{Ol, Ad, Gro}.

\par We investigate in this paper, alternatives for pseudofinite groups of the same flavour as the  Tits alternative. We show that an $\aleph_{0}$-saturated pseudofinite group either contains the free subsemigroup of rank $2$ or is supramenable. This follows from the following result: an $\aleph_{0}$-saturated pseudofinite group either contains the free subsemigroup of rank $2$ or is nilpotent-by-(uniformly locally finite) (Theorem 4.1). More generally we prove that an $\aleph_{0}$-saturated pseudofinite group satisfying a finite disjunction of Milnor identities is nilpotent-by-(uniformly locally finite) (Corollary 4.7). This is a straightforward consequence of the analogue proven in the class of finite groups (\cite{mil}).

\par Then we show that whether the following dichotomy holds for $\aleph_{0}$-saturated pseudofinite groups, namely it either contains a free nonabelian subgroup or it is amenable is equivalent to whether a finitely generated residually finite group which satisfies a nontrivial identity is amenable (respectively uniformly amenable) (Theorem \ref{prop-equiv}). 


 In the same spirit, we revisit  the results of S. Black \cite{B}  who considered a "finitary Tits alternative", i.e. an analog of Tits alternative for classes of finite groups. We reformulate  Black's results  in the context of pseudofinite groups  (Theorem \ref{thm-alter-bounded-rank}) and we strengthen it to the class of finite groups of weakly bounded rank.  A class of finite groups is {\it weakly of bounded rank} if the class of the radicals has bounded (Pr\"ufer) rank and the index of the sockels are bounded. We obtain the following dichotomies for an $\aleph_{0}$-saturated pseudo-(finite weakly of bounded rank) group $G$: either $G$ contains a nonabelian free group or $G$ is nilpotent-by-abelian-by-(uniformly locally finite) (Theorem 6.10). As S. Black, we use results of A. Shalev and D. Segal on classes of finite groups of bounded Pr\"ufer rank (\cite{Sh}, \cite{Se-rank}).
\par We will be also interested in classes of finite groups satisfying some uniform conditions on centralizer dimension, namely for which there is a bound on the chains of centralizers. A class $\C$ of finite groups has bounded $c$-dimension, if there is $d\in \N$ such that for each $G\in \C$ the $c$-dimension of $rad(G)$ and of the index of the sockels of $G/rad(G)$ are bounded by a function of $d$ only. We show that an $\aleph_{0}$-saturated pseudo-(finite of bounded $c$-dimension) group either contains a nonabelian free group or is soluble-by-(uniformly locally finite) (see Corollary 6.12). We use a result of E. Khukhro \cite{res} on classes of finite soluble groups of finite $c$-dimension. 
\par In our proofs, we use the following uniformity results which hold in the class of finite groups: the result of J. Wilson \cite{W}  who obtained a formula $\phi_{R}$ which defines across the class of finite groups the soluble radical, definability results for verbal subgroups of finite groups due to N. Nikolov and D. Segal \cite{Se} and the positive solution of the restricted Burnside problem due to E. Zelmanov \cite{V}.

The present paper is organized as follows. In the next section, we relate  the notion of being pseudofinite with other approximability properties by a class of groups and we recall some background material. In Section 3, we study some properties of finitely generated pseudofinite groups. Section 4 is devoted to the proof of the fact that $\aleph_{0}$-saturated pseudofinite group either contains the free subsemigroup of rank $2$ or is nilpotent-by-(uniformly locally finite) (Theorem 4.1).  Then, in Section 5 we study the general problem of the existence of nonabelian free subgroups and its relations with amenability. We end in Section 6 by giving the generalization (in the class of pseudofinite groups) of the above-mentionned Black's results and also some other alternatives under assumptions like bounded $c$-dimension.

\section{Generalities}
\par In this section we will first relate various notions of {\it approximability} of a group by a class of (finite) groups. The reader interested in a more thorough exposition can consult for instance the survey by T. Ceccherini-Silberstein and  M. Coornaert \cite{Coor-loc-sof}. We point out  that  Proposition \ref{prop-equiv-resid} (and its Corollaries)  seems new and it is important in the proof of Theorem \ref{prop-equiv}. At the end of this section we review some basic model-theoretic properties of pseudofinite groups.
\par In \cite{VG}, A.Vershik and E. Gordon considered a new version of embedding for groups; they defined $LEF$-groups, namely groups locally embeddable in a class of finite groups.  The definition adapts  to any class of groups  and it is related to various residual notions that we recall here.  
\nota Given a class $\C$ of $\L$-structures, we will denote by $Th(\C)$ (respectively by $Th_{\forall}(\C)$) the set of sentences (respectively
universal sentences) true in all elements of $\C$. 
\par Given a set $I$, an ultrafilter $\U$ over $I$ and a set of $\L$-structures $(C_i)_{i \in I}$, we denote by $\prod_I^\U C_i$ the ultraproduct of the family $(C_i)_{i \in I}$ relative to $\U$.  We denote by $\P_{fin}(I)$ the set of all finite subsets of $I$.  
\enota

\defn \label{def1} Let $\C$ be a class of groups.

$\bullet$ A group $G$ is called \textit{approximable} by  $\C$ (or locally $\C$ or locally embeddable into $\C$) if for any finite subset $F \subseteq G$, there exists a group $G_{F}\in \C$ and an {\it injective} map $\xi_{F}:F\rightarrow G_{F}$ such that $\forall g, h\in F$, if $gh\in F$, then $\xi_{F}(gh)=\xi_{F}(g)\xi_{F}(h)$. When $\C$ is a class of finite groups, then $G$ is called $LEF$.

$\bullet$ A group $G$ is called \textit{residually-$\C$}, if for any nontrivial element $g \in G$, there exists a homomorphism $\varphi : G \rightarrow C \in \C$ such that $\varphi(g) \neq 1$.

$\bullet$  A group $G$ is called \textit{fully residually}-$\C$, if for any finite subset  $S$ of nontrivial elements of $G$, there exists a homomorphism $\varphi : G \rightarrow C \in \C$ such that $1 \not \in \varphi(S)$.  

$\bullet$ A group $G$ is called \textit{pseudo-$\C$} if $G$ satisfies $Th(\C)=\bigcap_{C \in \C}Th(C)$.
\edefn
\par In particular, when $\C$ is the class of finite groups, a pseudo-$\C$ group is a pseudofinite group. In this case, we will abreviate pseudo-$\C$ group by pseudofinite group. We note that if $\C$ is closed under  direct product and subgroups, that is $\C$ is a pseudovariety,  then a group $G$ is residually-$\C$ if and only if it is fully residually-$\C$.

\par We will use the following variation of a theorem of Frayne  (\cite{CK} 4.3.13), which can be stated as follows. Let $\A$ be an $\L$-structure and $\C$ a class of $\L$-structures. Assume that $\A$ satisfies $Th(\C)$. Then there exists $I$ and an ultrafilter $\U$ on $I$ such that 
$\A$ elementarily embeds into an ultraproduct of elements of $\C$. 
It follows for instance that a group $G$ is pseudofinite if and only if it is elementarily embeddable in some ultraproduct of finite groups; a property that will be used throughout the paper without explicit reference.

\prop \label{prop1-frayne}
Let $\A$ be an $\L$-structure and $\C$ a class of $\L$-structures. Assume that $\A$ satisfies $Th_{\forall}(\C)$. 
Then there exists $I$ and an ultrafilter $\U$ on $I$ such that 
$\A$ embeds into an ultraproduct of elements of $\C$. 
\eprop
\pr We enumerate the elements of $A$ as $(a_{\alpha})_{\alpha<\delta}$ and we denote by $\L_{A}:=\L\cup\{c_{\alpha}:\;\alpha<\delta\}$. 
We will consider $\A$ as an $\L_{A}$-structure interpreting $c_{\alpha}$ by $a_{\alpha}$. Let $\Fa_{A}$ be the set of all $\L_{A}$-quantifier-free sentences $\phi(c_{\alpha_{1}},\cdots,c_{\alpha_{n}})$, where $\alpha_{1},\cdots,\alpha_{n}<\delta$.  Let $I:=\{\phi\in \Fa_{A}: \A\models \phi\}$. Note that if $\A\models \phi(c_{\alpha_{1}},\cdots,c_{\alpha_{n}})$, then there exists $\B\in \C$ such that $\B\models \exists x_{1}\cdots\exists x_{n}\;\phi(x_{1},\cdots,x_{n})$. Denote $\B_{\phi}$ such element of $\C$ and the corresponding tuple of elements $b_{\phi}:=(b_{\phi,\alpha_{1}},\cdots, b_{\phi,\alpha_{n}})$ such that $\B_{\phi}\models \phi(b_{\phi}).$
For any $\phi(c_{\alpha_{1}},\cdots,c_{\alpha_{n}})\in I$, we set $J_{\phi}:=\{\psi(c_{\alpha_{1}},\cdots,c_{\alpha_{n}})\in I:\;\B_{\psi}\models \phi(b_{\psi})\}$. These subsets $J_{\phi}$ have the finite intersection property and so there exists an ultrafilter $\U$ on $I$ containing these $J_{\phi}$.
\par Finally we define a map $f$ from $\A$ to $\prod_{I}^{\U} \B_{\phi}$ by sending $a_{\alpha}$ to $[b_{\phi\alpha}]_{\U}$ and check this is an embedding. Assume that for $\phi\in \Fa$, $\A\models \phi(a_{\alpha_{1}},\cdots,a_{\alpha_{n}})$, so $J_{\phi(c_{\alpha_{1}},\cdots,c_{\alpha_{n}})}\in \U$, so $\{\psi(c_{\alpha_{1}},\cdots,c_{\alpha_{n}})\in I:\B_{\psi}\models \phi(b_{\psi})\}\in \U$.
\qed
\medskip
\par One can derive from this proposition the following result of Malcev, namely that a group $G$ embeds in an ultraproduct of its finitely generated subgroups, by letting $\C$ to be the class of finitely generated subgroups of $G$. 
\prop \label{prop-equiva}Let $G$ be a group and $\C$ a class of groups. The following properties are equivalent.
\begin{enumerate}
\item The group $G$ is approximable by  $\C$.
\item $G$ embeds in an ultraproduct of elements of $\C$.
\item $G$ satisfies $Th_{\forall}(\C)$.
\item Every finitely generated subgroup of $G$ is approximable by $\C$. 

\end{enumerate}
\eprop
\pr  $(1) \Rightarrow (2)$.  Let $I=\P_{fin}(G)$ and let $\U$ be an  ultrafilter containing all subsets of the form $J_{F}:=\{e\in \P_{fin}(G): F\subset e\}$, with $F\in \P_{fin}(G)$.
 Choose $\xi_{F} : F \rightarrow G_F \in \C$ as in Definition \ref{def1}.  Then consider the ultraproduct $\prod^{\U}_{I} G_{F}$ and let for $g\in G$,  $\xi(g):=(\xi_{F}(g))_{F \in I}$. Then $\xi$ is a monomorphism.

\par $(2) \Rightarrow (1)$.  Assume that $G$ embeds in a  ultraproduct of elements of $\C$. Let $F\subset G$ a finite subset. We can describe the partial multiplication table of $F$ by a conjunction of basic formulas. Denote by $\si$ the existential sentence obtained by quantifying over the elements of $F$. This sentence is true on an infinite family of elements of $\C$. Let $H\in \C$ satisfying this sentence and define a map from $G$ to $H$ accordingly.

\par $(2) \Rightarrow (3)$.  Let $\si\in Th_{\forall}(\C)$. Then since $G$ embeds in an ultraproduct of elements of $\C$, $G\models \si$.

The implication  $(3) \Rightarrow (2)$ is the statement of   Proposition \ref{prop1-frayne} and we see also that the equivalence $(1) \Leftrightarrow (4)$ is clear.  \qed

\prop \label{prop-equiv-resid} Let $G$ be a group and $\C$ a class of groups. The following properties are equivalent.
\begin{enumerate}
\item The group $G$ is approximable by  $\C$.
\item  For every finitely generated subgroup $L$ of $G$, there exists a sequence of finitely generated residually-$\C$ groups $(L_n)_{n\in \mathbb N}$ and a sequence of homomorphisms $(\varphi_n : L_n \rightarrow L_{n+1})_{n\in \mathbb N}$ such the following properties holds:

\smallskip
$(i)$  $L$ is the direct limit, $L= \varinjlim L_{n}$, of the system $\varphi_{n,m} : L_n \rightarrow L_m$, $m \geq n$,  where $\varphi_{n,m}=\varphi_m \circ \varphi_{m-1} \cdots \circ \varphi_n$. 

$(ii)$ For any $n \geq 0$, for any finite subset $S$ of $L_n$, if $1 \not \in \psi_n(S)$, where $\psi _n : L_n \rightarrow L$ is the natural map, there exists a homomorphism $\varphi : L_n \rightarrow C \in \C$ such that $1 \not \in \varphi(S)$.

\end{enumerate}
\eprop

\pr  $(1) \Rightarrow (2)$.   Let $L$ be a finitely generated subgroup of $G$.  Let  $$L=\<a_1, \dots, a_p\,| r_0, \dots, r_n, \dots\>$$ be a presentation of $L$ and set $D_n=\<a_1, \dots, a_p\,| r_0, \dots,r_n\>$ for $n \geq 0$. Let $\<x_{1},\cdots,x_{p}|\>$ be the free group generated by $x_{1},\cdots,x_{p}$ and $N_{n}$ the normal subgroup generated by $r_0(\bar x), \dots, r_n(\bar x)$, where $\bar x:=(x_{1},\cdots,x_{p})$. Then $D_{n}\cong \<x_1, \dots, x_p\>/N_{n}$.  We have a direct system of homomorphisms $f_{n,m}$ from $D_{n}$ to $D_{m}$, $n\leq m$, defined by $f_{n,m}(x.N_{n})=x.N_{m}$. It follows that  $L$ is the direct limit of the previous  system, $L= \varinjlim D_{n}$.

Let us define   $L_n$ to be the group $D_n/K_{n}$,  where $K_{n}$ is the intersection of all normal subgroups $M$ for which $D_n/M$ is a subgroup of some $C \in \C$.  We see that  each $L_n$ is residually-$\C$.  Let $\pi_n : D_n \rightarrow L_n$ be the natural homomorphism. 

Clearly, we have a natural homomorphism $\varphi_{n,m} : L_n \rightarrow L_{m}$ such that $\pi_m \circ f_{n,m}=\varphi_{n,m} \circ \pi_{n}$. We let $\C L$ to be the direct limit of the given system,  $\C L=\varinjlim L_n$. We note also that, we have a natural homomorphism $\pi : L \rightarrow \C L$ and we get the following diagram.

$$\xymatrix{
    D_0 \ar[r]^{f_{0,1}} \ar[d]_{\pi_0}  & L_1 \ar[d]^{\pi_1}  \ar[r]& \cdots \ar[r] &D_n \ar[d]^{\pi_n} \ar[r]^{f_{n, n+1}}&D_{n+1} \ar[r] \ar[d]^{\pi_{n+1}}&\cdots& L \ar[d]^\pi\\
     L_0 \ar[r]_{\varphi_{0,1}} & \ar[r] L_1&\cdots \ar[r]&L_n \ar[r]_{\varphi_{n, n+1}}&L_{n+1} \ar[r]&\cdots& \C L
  }$$

 We claim that $\pi$ is an isomorphism. By definition $\pi$ is surjective and it is sufficient to show that it is injective.  We note that for any word $w(\bar a)$, $\pi(w(\bar a))=1$ if and only if there exists $n \in \N$ such that $\pi_n(w(\bar a))=1$.

Let $a\in L\setminus\{1\}$. Then there is a word $w$ in $\bar a$ such that $a=w(\bar a)$ and so for all $m\in \N$, $L\models \exists \bar x\,(w(\bar x)\neq 1\;\&\;\bigwedge_{0\leq n\leq m} r_{n}(\bar x)=1)$.

By hypothesis, $L$ is approximable by $\C$, so for all $m\in \N$, there exists $C_{m}\in \C$ such that $C_{m}\models \exists \bar x\,(w(\bar x)\neq 1\;\&\;\bigwedge_{0\leq n\leq m} r_{n}(\bar x)=1)$.
Let $\bar b_{m}\in C_{m}$ such that $w(\bar b_{m})\neq 1$ and $\bigwedge_{0\leq n\leq m} r_{n}(\bar b_{m})=1$.
Hence, there is a homomorphism  from $L_{m}$ to the subgroup of $C_{m}$ generated by $\bar b_m$, so for some normal subgroup $M_{m}$ of $L_{m}$, we get $L_{m}/M_m\cong \<\bar b_{m}\> \leq C_m$.  By definition, we have $K_m \leq M_m$ and thus $\pi_m(w(\bar a))\neq 1$  (for all $m \in \N)$. Hence $\pi(a) \neq 1$ and thus $\pi$ is injective as required.

Let $n \geq 0$ and $ S=\{g_1, \dots, g_q\}\subseteq L_n$ be a finite subset such that $1 \not \in \psi_n(S)$, where $\psi _n : L_n \rightarrow \C L$ is the natural map.  Proceeding as above, there exists a finite sequence of words $(w_j(\bar x))_{1 \leq j \leq q}$ such that $g_j=w_j(\bar x)$ and $L\models \exists \bar x\,(\bigwedge_{1 \leq j \leq q}w_j(\bar x)\neq 1\;\&\;\bigwedge_{0\leq n\leq m} r_{n}(\bar x)=1)$.  Procceding as above, we find a normal subgroup $M_n\leq L_n$ such that $K_n  \leq M_n$, $g_j \not \in M_n$,  $L_n/M_n$ isomorphic to a subgroup of some element $C \in \C$; which gives the required result.

\par $(2) \Rightarrow (1)$.   Let $L$ be a finitely generated subgroup of $G$. Let $(L_n)_{n \in \Bbb N}$ be a sequence of  residually-$\mathcal C$ groups whose direct limit is $L$ and satisfying the property $(ii)$.  Denote the maps in the direct system between $L_{n}$ and $L_{m}$, $n\leq m$, by $f_{n,m}$.  Let $L=\varinjlim L_{n}=\bigsqcup_{n} L_{n}/\sim$, and for  $x\in L$ we let  $x_{n}\in L_{n}$ to be  a representative of $x$ with respect to the equivalence relation $\sim$.

Then we see that if $L$ satisfies a formula of the  form, for some tuple $\bar x$, 
$$
\bigwedge_{1 \leq j \leq q}w_j(\bar x)\neq 1\;\&\;\bigwedge_{1\leq i\leq p} r_{i}(\bar x)=1,
$$
then there exists $n \in \mathbb N$ such that 
$$
L_n \models  \bigwedge_{1 \leq j \leq q}w_j(\bar x_n)\neq 1\;\&\;\bigwedge_{1\leq i\leq p} r_{i}(\bar x_n)=1, 
$$
where $\bar x_n$ is a representative of $\bar x$. By $(ii)$, we conclude that there exists $C_n \in \C$ such that 
$$
C_n  \models \exists  \bar x \; (\bigwedge_{1 \leq j \leq q}w_j(\bar x)\neq 1\;\&\;\bigwedge_{1\leq i\leq p} r_{i}(\bar x)=1). 
$$

We conclude that $L$ satisfies $Th_{\forall}(\C)$. By the result of Malcev recalled above, $G$ embeds in an ultraproduct of its finitely generated subgroups and so $G\models Th_{\forall}(\C)$. \qed
\medskip

\kor  Let $G$ be a group and $\C$ a pseudovariety  of groups. Then $G$ is approximable by $\C$ if and only if any finitely generated subgroup of $G$ is a direct limit of finitely generated fully residually-$\C$ groups. \qed

\ekor

As consequence, taking $\C$ to be the class of finite groups, we have the following corollary which  seems new and not observed in the literature. A. Vershik and E. Gordon showed that a finitely presented LEF-group is residually finite \cite{VG}.

\kor \label{cor1} Let $G$ be a finitely generated group. The following properties are equivalent. 

\par $(1)$   $G$ is LEF. 

\par $(2)$  $G$ is a direct limit of residually finite groups. \qed

\ekor

\bem $\;$

$(1)$  We note also  that a finitely generated group is approximable by $\C$ if and only if $G$ is a limit in an adequate topological space of marked groups (see \cite{Coor-loc-sof}, \cite{CG}). 

$(2)$  It follows from Proposition \ref{prop-equiva}  that the class of pseudofinite groups is included into the class of $LEF$-groups since any pseudofinite group embeds into an ultraproduct of finite groups. It is easy to see that the class of pseudofinite groups is strictly smaller than the class of $LEF$-groups (see below).
\ebem

\noindent {\bf Examples.} $\;$
\begin{enumerate}
\item Let $\C$ be the class of finite groups.  A locally residually finite group is locally $\C$ \cite{VG}. There are groups which are not residually finite and which are approximable by  $\C$,  for instance,  in \cite{Coor-loc-sof} an example of a finitely generated amenable LEF group which is not residually finite is given. There are residually finite groups which are not pseudofinite, for instance the free group $F_{2}$ (see Corollary \ref{trf-hyp-pseudo}).
\item Let $\C$ be the class of free groups. 
If $G$ is  fully residually-$\C$ (or equivalently $\omega$-residually free or a limit group), then $G$ is approximable by $\C$ \cite{Chis}. Conversely if $G$ is approximable by $\C$, then $G$ is locally fully residually-$\C$. The same property holds also in hyperbolic groups \cite{sela-hyp, weidmann-equa} and more generally in equationally noetherian groups \cite{Ould-equa}. 
\item Let $V$ be a possibly infinite-dimensional vector space over a field $K$. Denote by $GL(V,K)$ the group of automorphisms of $V$. Let $g\in GL(V,K)$, then $g$ has finite residue if the subspace $C_{V}(g):=\{v\in V: g.v=v\}$ has finite-co-dimension. A subgroup $G$ of $GL(V,K)$ is called a {\it finitary} (infinite-dimensional) linear group, if all its elements have finite residue. A subgroup $G$ of $\prod_{i\in I}^{\U} GL(n_{i},K_{i})$, where $K_{i}$ is a field, is {\it of bounded residue} if for all $g\in G$, where $g:=[g_{i}]_{\U}$, $res(g):=inf \{n\in \N: \{i\in I: res(g_{i})\leq n\}\in \U\}$ is finite.
 \par E. Zakhryamin has shown that any finitary (infinite-dimensional) linear group $G$ is isomorphic to a subgroup of bounded residue of some ultraproduct of finite linear groups (\cite{Zak} Theorem 3).
In particular letting $\C:=\{GL(n,k)$, where $k$ is a finite field and $n\in \N\,\;\}$, any finitary (infinite-dimensional) linear group $G$ is  approximable by $\C$.
\end{enumerate}
\bigskip
\par Recall that a group $G$ is said to be a \textit{CSA-group} \cite{MR} if every maximal abelian subgroup $A$ of $G$ is malnormal; that is $A^g\cap A=1$ for any $g \in G\setminus A$. In particular, a nonabelian CSA-group has no nontrivial normal proper abelian subgroup. Let us observe that if $G$ is $CSA$, then all the centralizers are abelian. Indeed,  let $a\in G\setminus\{1\}$ and  let $A$ be  a maximal abelian subgroup of $G$ containing $a$ and suppose that there exists $b\in C_{G}(a)\setminus A$. Consider $A^b\cap A$. This intersection contains $a$, which is a contradiction. In particular the maximal abelian subgroups of $G$ are centralizers.
\lmm {\rm \cite{MR, ould-csa}} \label{CSA} The property of being CSA can be expressed by a universal sentence.
\elmm
\pr Let us express that $\forall x\neq 1\;C(x)$ is abelian and $\forall y \forall z\; y\notin C(x)$ and $z\in C(x)\cap C(x)^y$ implies that $z=1$. Then $G$ is $CSA$ iff $G$ satisfies that sentence.
\qed
\kor {\rm  \cite{ould-csa}} A finite CSA group is abelian.
\ekor
\pr Since the property of being CSA is universal, it is inherited by subgroups. So, a minimal nonabelian CSA finite group has all its proper subgroups abelian and so this group is soluble by a result of O.J.Smidt (\cite{R} 9.1.9) and thus it has a nontrivial proper normal abelian subgroup; a  contradiction.
\qed
\prop \label{CSA-abelian} A  pseudofinite CSA-group is abelian.  \eprop

\pr Indeed, $G \preceq \prod^{\U}_I F_{i}$, where each $F_i$ is finite,  and since the class of CSA-groups is axiomatizable by a single universal sentence (Lemma \ref{CSA}), for almost $i$, $F_i$ is a CSA-group. But a finite CSA-group is abelian and thus $G$ is abelian.   \qed
\kor The classes of pseudofinite groups and of nonabelian groups approximable by nonabelian free groups have a trivial intersection.\qed
\ekor
\pr A nonabelian free group is a CSA-group and we apply Lemma \ref{CSA-abelian} and Proposition \ref{prop-equiva}.
\qed
\bigskip

\par There are  other  kinds  of approximation by classes of groups related  to the previous notions.  Gromov \cite[Sect.6.E]{G} introduced groups whose Cayley graphs are  {\it initially subamenable} which are afterward called  {\it sofic} by B. Weiss \cite{weiss}. These groups can  be seen as a simultaneous generalization of amenable groups and residually finite groups. We give a definition which is a slight generalization of already known notions by using {\it invariant metric} and which follows the definition given in \cite{Lev-Louis} (see also \cite{Lev-Louis} for the proof of the fact that this definition is equivalent to the classical one for sofic groups).   A group $G$ is an \textit{invariant-metric group}, if there is a distance $d$ on $G$     which is  bi-invariant; namely for any $x,y,  z \in G$, $d(zx, zy)=d(xz, yz)=d(x,y)$. 

\defn \label{def2} Let $\C$ be a class of invariant-metric groups.  A group $G$ is {\it $\C$-sofic or sofic relative to $\C$},  if for any finite subset $F$ of $G$, there exists $\epsilon >0$ such that for every $n\in \N^*$, there exists $(C, d_C) \in \C$ and an injective map $\xi_{F}:F\rightarrow C$ such that  for any $g,h \in F$, if $gh\in F$, then $d_{C}(\xi_{F}(gh),\xi_{F}(g)\xi_{F}(h))\leq \frac{1}{n}$ and for all $g\in F\setminus \{1\}$,  $d_{C}(1,\xi_{F}(g))\geq \epsilon$. 
\edefn

For $n \in \N^*$ let $S_n$ be the symmetric group on $n$ elements and $d_n$ be the distance on $S_n$, called  the {\it normalized Hamming distance},  defined by $d_{n}(\si,\tau)=\frac{1}{n}\vert \{i\in n: \si(i)\neq \tau(i)\}\vert$ with $\si,\;\tau\in S_n$ (identifying $n$ with the subset $\{1,\cdots,n\}$ of natural numbers).  Then  a sofic group relative to $\C=\{(S_n, d_n):\,n \in \N\}$ is called \textit{sofic}. 

We are interested in a characterization of sofic groups relative to $\C$  in terms of embeddings in adequate ultraproducts analogous to that   of  Proposition \ref{prop-equiva}. Elek-Szab\'o \cite{Elek-Saz}  gave   such characterization for  sofic groups  that we generalize  here to the  general framework of  invariant-metric  groups (see also \cite{Pe, weiss, Coor-loc-sof}).

\defn \label{def3} Let $\C$ be a class of invariant-metric groups, $I$ a set and $\U$ a nonprincipal ultrafilter on $I$.  For a sequence $(C_i)_{i\in I}$ from $\C$ we let $\mathcal G=\prod_I^\U C_i$.  Then $\mathcal G$ is group endowed with  a natural bi-invariant metric $d_\U$ with values in $\prod_I^\U \R$ by defining $d_\U([a_i]_\U, [b_i]_\U)=[d_{C_i}(a_i, b_i)]_\U$.  Consider the subset of $\mathcal G$ defined by $\mathcal N=\{g \in  \mathcal G| d_\U(1, g) \hbox{ is infinitesimal} \}$. Then $\mathcal N$ is a normal subgroup and the quotient group $\mathcal G/ \mathcal N$ will be called an {\it universal $\C$-sofic group}.  
\edefn

\prop Let $\C$ be a class of invariant-metric groups and $G$ a group. Then the following properties are equivalent. 

$(1)$ $G$ is $\C$-sofic. 

$(2)$ $G$ is embeddable in some  universal $\C$-sofic group. 

\eprop

\proof

$(1)  \Rightarrow (2)$.  Let $J(G)=\P_{fin}(G) \times \Bbb N$ and let $\U$ be a nonprincipal ultrafilter over  $J(G)$ containing all subsets of the form $J_{F, n_0}=\{(e,n)\in J:\;F\subset e\;\&\;n\geq n_{0}\}$. For each $(e,n) \in J(G)$ let $C_{(e,n)} \in \C$ and $\xi_{e,n} : e  \rightarrow C_{(e,n) }$ for which the properties given in Definition \ref{def2} are fulfilled.  Consider the ultraproduct $\mathcal G(G)=\prod_{J(G)}^\U C_{(e,n)}$ and $\mathcal N(G)$ the corresponding normal subgroup as defined in Definition \ref{def3}.   Define $\xi : G \rightarrow \mathcal G(G)$ by $\xi(g)=[\xi_{(e,n)}(g))_{(e,n) \in J(G)}]_\U$.   We note that $d_\U(\xi(g_1g_2), \xi(g_1)\xi(g_2))$ is infinitesimal and $d_{\U}(\xi(g),1) >0$ for every $g \in G \setminus \{1\}$.  Hence $\xi : G \rightarrow  \mathcal G(G)/ \mathcal N(G)$ is an embedding.

$(2) \Rightarrow (1)$.  Let $\xi : G \rightarrow  \mathcal G/ \mathcal N$ be an embedding and set $\pi : \mathcal G \rightarrow \mathcal G/ \mathcal N$ the natural map. For every finite set $F$ of $G$, let $F'$ be a subset of $\mathcal G$ such that the restriction of $\pi$ to $F'$ is a bijection from $F'$ to $\xi(F)$ and set $\pi_F^{-1} :\xi(F) \rightarrow F'$ the inverse of the restriction of $\pi$. Then  for any $g, h \in F$, if $gh \in F$, then $d_\U(\pi_F^{-1}\circ \xi (gh), \pi_F^{-1}\circ \xi(g) \cdot \pi_F^{-1} \circ \xi(h))$ is infinitesimal and $d_\U(1, \pi_F^{-1} \circ \xi(f))\geq \epsilon$ for any $f  \in F$ and some  $\epsilon>0$. By considering an adequate subset of $\U$, we get the required conclusion.
\qed

\bigskip

It is an open problem  whether or not any group is sofic. However, it is known that many groups are sofic: residually finite  groups, LEF-groups, amenable groups,  residually amenable groups (see for instance \cite{Coor-loc-sof}). More generally any pseudosofic group is sofic \cite[Proposition 7.5.10]{Coor-loc-sof}. We see in particular that any pseudofinite group is sofic.

\bigskip
  
\par The next lemma is well-known and holds for any pseudofinite structure, but we give a proof for the reader convenience.
\lmm \label{def} Let $G$ be a pseudofinite group. Any definable subgroup or any
quotient by a definable normal subgroup is pseudofinite. 
\elmm
\pr Since $G$ is pseudofinite, there is a family $(F_i)_{i \in I}$ of finite groups and an ultrafilter $\U$ such that $G \preceq \prod^{\U}_I F_{i}$.  Let $\phi(x, \bar y)$ be a formula and $\bar b \in G$ such that $\phi(x,\bar b)$ defines a subgroup of $G$. Let $[\bar b_i]_{i \in I}$ be a sequence representing  $\bar b$.  Then there exists $U \in \U$ such that for every $i \in U$, $\phi(x, \bar b_i)$ defines a subgroup of $F_i$. 

Given  a formula  $\theta(\bar x)$ whose prenex form is: $Q z_{1}\cdots Q z_{n}\;\chi(\bar z,\bar x)$, where $\chi$ is a quantifier-free formula and $Q$ denotes either $\exists$ or $\forall$, we let   $\theta^{\phi}(\bar x;\bar y)=Q z_{1}\cdots Q z_{n}\;\chi(\bar z,\bar x)\;\&\;\bigwedge_{i=1}^n \phi(z_{i},\bar y).$

\par Assume now that $\sigma$ is a sentence. Then $\sigma^{\phi}(\bar y)$ expresses that the subgroup defined by $\phi(x; \bar y)$ satisfies $\sigma$. We conclude that $\prod^{\U}_I F_{i} \models \sigma^{\phi}(\bar b)$ for any sentence $\sigma$ true in  any finite group. Hence $\phi(G; \bar b)$ satisfies any sentence true in any finite group, and thus it is pseudofinite.  The same method works for quotients by using relativization to quotients.  This time instead of relativizing the quantifiers, we have to replace equality by belonging to the same coset of $\phi(G;\bar b)$.
\qed

\bigskip
\par There are many definitions of semi-simple groups in the literature, which differ from a context to another. We adopt here  the following. We will say that a group is \textit{semi-simple} if it has no nontrivial normal abelian subgroups. We note in particular that a semi-simple group has no nontrivial soluble normal subgroups. 

 \par  Let $G$ be a finite group and let $rad(G)$ be the soluble radical, that is the largest normal soluble subgroup of $G$. In \cite{W}, J. Wilson proved the existence of a formula, that will be denoted in the rest of this paper by $\phi_{R}(x)$,  such  that  in any finite group $G$, $rad(G)$ is definable by $\phi_R$. 
\lmm \label{radical} If  $G$ is a pseudofinite group then $G/\phi_{R}(G)$ is a pseudofinite semi-simple group.
\elmm
\pr By the preceding Lemma \ref{def} and the above result of Wilson, $G/\phi_{R}(G)$ is pseudofinite. Let us show it is semi-simple. Suppose there exists $a\in G\setminus \phi_{R}(G)$ such that for all $h,\;g \in G$,  $\phi_{R}([a^h,a^g])$. By hypothesis $G\equiv \prod^{\U}_I G_{i}$, where each $G_i$ is finite. So on an element of $\U$, $G_{i}\models \exists x\forall y\forall z\;\phi_{R}([x^y,x^z])\;\&\;\neg\phi_{R}(x)$; a contradiction.
\qed
\medskip

\par We end the section by recalling another well-known result, namely the equivalence in an $\aleph_{0}$-saturated group of not containing the free group and of satisfying a nontrivial identity.
\nota  Let $F_{2}$ be  the free nonabelian group on two generators and  $M_{2}$ be the free subsemigroup on two generators.
\enota
\lmm \label{lem-dicho}Let $G$ be an $\aleph_{0}$-saturated group. Then either $G$ contains $F_{2}$, or $G$ satisfies a nontrivial identity (in two variables). In the last case, either $G$ contains $M_{2}$, or $G$ satisfies a finite disjunction of positive nontrivial identities in two variables.
\elmm
\pr We enumerate the set $W_{x,y}$ (respectively $M_{x,y}$) of nontrivial reduced words in $\{x, y,x^{-1},y^{-1}\}$ (respectively in $\{x, y\}$) and we consider the set of atomic formulas $p(x,y):=\{\theta(x,y):=(t(x,y)\neq 1): t(x,y)\in W_{x,y}\}$ (respectively $q(x,y):=\{\theta(x,y):=(t_{1}(x,y)\neq t_{2}(x,y)): t_{1}(x,y)\,,t_{2}(x,y)\in M_{x,y}\cup\{1\}, t_{1}\neq t_{2}\}$).
\par Either there is a finite subset $I$ of $p(x,y)$ (respectively of $q(x,y)$) which is not satisfiable in $G$ and so $G\models \forall x\;\forall y\;\bigvee_{\theta\in I} \theta(x,y)$, otherwise since $G$ is $\aleph_{0}$-saturated, $G\supset F_{2}$ (respectively $G\supset M_{2}$).
\par Observe that if a group $G$ satisfies a finite disjunction of nontrivial identities in two variables, then it satisfies one nontrivial identity.
For sake of completeness, let us recall here the argument. Suppose $G\models (t_{1}(x,y)=1\vee t_{2}(x,y)=1)$. Either $t_{1}(x,y)$ and $t_{2}(x,y)$ do not commute in the free group generated by $x,\;y$ and so the commutator $[t_{1},t_{2}]\neq 1$ in the free group and so the corresponding reduced word is nontrivial and $G\models [t_{1},t_{2}]=1$. Or $t_{1},\;t_{2}$ do commute in the free group and so there exists a nontrivial reduced word $t$ in $x,\;y$ and $z_{1},\;z_{2}\in \Z$ such that $t_{1}=t^{z_{1}}$ and $t_{2}=t^{z_{2}}$. In that last case $G\models t(x,y)^{z_{1}.z_{2}}=1$.
\qed

\section{Finitely generated pseudofinite groups.}

We study in this section some properties of finitely generated pseudofinite groups, led by a question of G. Sabbagh who asked whether all such groups were finite. There will be two main ingredients. First, a definability result due to N. Nikolov and D. Segal that we will recall below (Theorem \ref{Segal}), and the following observation.
\par Recall that a group $G$ is said to be \textit{Hopfian} if any surjective endomorphism of $G$  is bijective. Any finitely generated residually finite group is Hopfian (Malcev) (see for instance \cite{MKS} page 415).  Since in a finite structure, any injective map is surjective and 
vice-versa, any definable map (with parameters) in a pseudofinite
group is injective iff it is surjective. In particular a pseudofinite group is \textit{definably} hopfian; that is any definable injective homomorphism is surjective. 
\nota Let $G^n$ be the verbal subgroup of $G$ generated by the set of all
$g^n$ with $g\in G$, $n\in \N$. The width of this subgroup is the maximal number (if finite) of $n^{th}$-powers necessary to write an element of $G^n$.
\enota
\defn A group $G$ involves a subgroup $H$, or $H$ is a section of $G$,  if there are subgroups $B \leq A\leq G$ with $B$  normal in  $A$ such that $A/B\cong H$.
\edefn
\thm {\rm \cite{Niko}, \cite{Se}(4.7.5.)\label{Segal} }There exists  a function $d  \rightarrow c(d)$,  such that if  $G$ is a $d$-generated finite group  and $H$ is a normal subgroup of $G$, then every element of  $[G, H]$ is a product of at most  $c(d)$ commutators of the form  $[h,g ]$,  $h \in H$ and  $g \in G$.  \par Moreover, in a finite group generated by $d$ elements not involving the alternating groups $A_{m}$ for $m\geq s$, the verbal subgroup generated by the $n^{th}$-powers, is of finite width bounded by $b(s,d,n)$.  \qed
\ethm
\par Finally let us recall the positive solution of the restricted Burnside problem, a long standing problem that was completely solved by E. Zemanov (\cite{V}, \cite{Zel-bur}). Given $k, d$, there are only finitely many finite groups generated by $k$ elements of exponent $d$.
\medskip
\prop \label{sab}{\rm (Sabbagh)} Any abelian  finitely generated pseudofinite group is finite. 
\eprop
\pr A finitely generated abelian group is a direct sum of a finite group and finitely many copies of $\Z$. 
So there exists a natural number $n$ such that $G^n$ is a $0$-definable subgroup of $G$ which is isomorphic to $\Z^k$ for some $k$.
But $\Z^k$ cannot be pseudofinite since the map $x\rightarrow x^2$ is injective
but not surjective.
\qed
 
\kor \label{trf-hyp-pseudo}There are no nontrivial torsion-free hyperbolic pseudofinite groups. \ekor
\pr A torsion-free hyperbolic group is a CSA-group and thus if it were pseudofinite then it would be abelian by Proposition \ref{CSA-abelian}.  We conclude by the above proposition.   \qed

\bigskip

\lmm \label{A_{n}} Suppose that there exists an infinite set $U \subseteq \Bbb N$  such that for any $n \in U$, the  finite group $G_{n}$ involves $A_{n}$. Then for any non-principal ultrafilter $\U$ containing $U$,  $G:=\prod^{\U}_{\N} G_{n}$ contains $F_{2}$.
\elmm
\pr Since $G$ is $\aleph_0$-saturated, using Lemma 2.7, it suffices to show that such group does not satisfy any nontrivial identity. By the way of contradiction let $w(x,y)=1$ be a nontrivial identity satisfied by $G$. Then letting $S_{n}$ be the full permutation group on $n$ letters, we would have that for infinitely many $n$, $S_{n}$ would satisfy the identity $w(x^2,y^2)=1$.  But then any finite group would satisfy a nontrivial identity since it embeds in some $S_{n}$ for $n$ sufficiently large. However $F_{2}$ is residually finite and so $F_{2}$ would satisfy a nontrivial identity, a contradiction.
\qed
\bigskip
\par Recall that a group is said to be \textit{uniformly locally finite} if for any $n \geq 0$,  there exists $\alpha(n)$ such that any $n$-generated subgroup of $G$ has cardinality bounded by $\alpha(n)$. In particular an uniformly locally finite group is of finite exponent.  Examples of uniformly locally finite groups include $\aleph_0$-categorical groups.

\lmm \label{res-burn} A pseudofinite group of finite exponent is uniformly locally finite.
\elmm
\pr Let $\<g_{1},\cdots,g_{k}\>$ be a $k$-generated subgroup of $G$. By definition $G \preceq  \prod^{\U}_J G_{j}$, where $G_{j}$ is a finite group. If $G$ is of exponent $e$, on an element of $\U$, $G_{j}$ is of exponent $e$. Let $[g_{mj}]_{j\in J}$, $1\leq m\leq k$, be a representative for $g_{m}$ and consider the subgroup $\<g_{1j},\cdots,g_{kj}\>$ on that element of $\U$. Then by the positive solution of the restricted Burnside problem, there is a bound $N(k,e)$ on the cardinality of that subgroup.
So the subgroup $\<g_{1},\cdots,g_{k}\>$ embeds into an ultraproduct of groups of cardinality bounded by $N(k,e)$ and so has cardinality bounded by $N(k,e)$.
\qed
\kor\label{app-burn} A group $G$ approximable by a class $\C$ of finite groups of bounded exponent is uniformly locally finite.
\ekor
\pr By Proposition \ref{prop-equiva}, such group embeds in an ultraproduct of elements of $\C$. So by the same reasoning as in the above lemma, any subgroup of $G$ generated by $k$ elements embeds into an ultraproduct of groups of cardinality bounded by a natural number $N(k,e)$ where $e$ is a bound on the exponent of the elements of $\C$ and so it is finite. \qed
\prop \label{prop-def} Let $L$ be a pseudo-($d$-generated finite groups). Then for any definable subgroup $H$ of $L$, the subgroup $[H, L]$ is definable. In particular for any $n \geq 0$, the derived subgroup $L^{(n)}$ is $0$-definable and of  finite width.  Similarly for terms of the descending central series of $L$. If moreover $L$ is $\aleph_0$-saturated and  does not contain $F_{2}$, then the verbal subgroups $L^n$, $n\in \N^*$, are $0$-definable  of finite width and of finite index.
\eprop

\pr Let $L\preceq \prod_I^{\U} F_{i}$, where $F_{i}$ is a finite group generated by $d$ elements. 
\par Let $\phi(x; \bar y)$ be a formula and $\bar b=[\bar b_i]$ such that $\phi(x; \bar b)$ defines a subgroup $H$.  On an element of the ultrafilter, $\phi(x; \bar b_i)$ defines a subgroup $H_i$ and the  subgroup $[H_i,F_{i}]$ is of width $\leq c(d)$ (Theorem \ref{Segal}). This property can be expressed by a sentence 
$$\bigwedge_{1\leq j \leq c(d)+1}\forall u_{j}\forall v_{j} \bigwedge_{1 \leq j \leq c(d) }\exists x_j \exists y_j(\bigwedge_{1 \leq j \leq c(d)} \phi(x_j;\bar b_i) \& \bigwedge_{1 \leq j \leq c(d)+1} \phi(u_j;\bar b_i) \Rightarrow \prod_{j=1}^{c(d)+1} [u_j,v_j]=\prod_{i=1}^{c(d)} [x_{j},y_{j}]),$$
and so, $[H_{i},F_{i}]$ is definable
as well as   the  subgroup $[H, L]$ of $L$.  A similar argument shows that  any term of the derived series of $L$ is definable and of finite width, as well as terms of the descending central series.

Suppose that $L$ is $\aleph_0$-saturated and  does not contain $F_{2}$. Then $L$ satisfies a nontrivial identity (Lemma \ref{lem-dicho}), as well as $\prod_I^{\U} F_{i}$. Hence by Lemma \ref{A_{n}}, there exists $s$ such that on an element of the ultrafilter $\U$,  $F_{j}$ does not involve $A_{m}$ for any $m\geq s$.   Therefore by Theorem \ref{Segal}, the sentence $\forall u \forall u_{1}\cdots\forall u_{b(s,d,n)}\;\exists x_{1}\cdots\exists x_{b(s,d,n)}\;\,u^n.\prod_{i=1}^{b(s,d,n)} u_{i}^n=\prod_{i=1}^{b(s,d,n)} x_{i}^n$ holds in $ \prod_I^{\U} F_{j}$. Since it holds in $\prod_I^{\U} F_{j}$, it holds in $L$ and so $L^n$ is $0$-definable and of finite width. 
\par By the solution of the restricted Burnside problem, the index of $F_j^n$ in $F_j$ is bounded in terms of $d$ and $n$ only. Then one can express that property by a $\exists\forall\exists$-sentence which  transfers in the ultraproduct of the $F_{j}$'s and therefore in $L$.
\qed

\nota Let $G$ be a group and $a, b\in G$, let $n\in \N$. We denote by $B^G_{\{a,b,a^{-1},b^{-1}\}}(n)$ the set of elements of $G$ which can be written as a word in $a,b,a^{-1},b^{-1}$ of length less than or equal to $n$. By convention the identity of the group is represented by a word of length $0$.
\enota
\defn \cite{B} A (finite) group $G$ contains an approximation of degree $n$ to $F_{2}$, the free nonabelian group on two  generators $x,y$ if there exists $a,b\in G$ such that $\vert B^G_{\{a,b,a^{-1},b^{-1}\}}(n) \vert=\vert B^{F_{2}}_{\{x, y,x^{-1},y^{-1}\}}(n)\vert$.
\edefn
\nota Let $G,\; L$ be two groups. Then $G\preceq_{\exists} L$ if $G$ is a subgroup of $L$ and every existential formula with parameters in $G$ which holds in $L$, holds in $G$.
\enota
\prop \label{prop-fg-derived}Let $G$ be a finitely generated pseudofinite group. Then the terms of the derived series are  $0$-definable of finite width and of finite index. 
If in addition $G$ does not contain any approximation of degree $n$ to $F_{2}$, $n\in \N$, then the sugbroups $G^m$ are $0$-definable of finite width and of finite index, $m\in \N^*$.
\eprop

\proof  Since $G$ is pseudofinite, $G \preceq L=\prod_I^\U G_i$, where each $G_i$ is finite. Let $\bar a$ be a finite generating tuple of $G$ and set  $\bar a=[\bar a_i]$,  $F_i=\<\bar a_i\>$ the subgroup of $G_i$ generated by $\bar a_i$.  We see  that $G \preceq_{\exists} \prod_I^\U F_i$.  Since $\prod_I^\U F_i$ is pseudo-($d$-generated finite groups), as in the proof of Proposition \ref{prop-def} any element of the derived subgroup is a product of at most  $c(d)$ commutators. Since this can be expressed by $\forall \exists$-sentence and as $G \preceq_{\exists} \prod_I^\U F_i$, we conclude that the same property holds in $G$, and thus $[G, G]$ is $0$-definable and of finite width.
\par By Lemma  \ref{def}, $G/[G,G]$ is a finitely generated pseudofinite abelian group, and so by Proposition \ref{sab}, it is finite. Hence $[G, G]$ is finitely generated and since it is $0$-definable, it is again pseudofinite (Lemma \ref{def}). Thus the conclusion on the terms of the derived series follows by induction. 

If $G$ does not contain any approximation of degree $n$ to $F_{2}$ then $\prod_I^\U F_i$ does not contain a free nonabelian group (Lemma \ref{lem-dicho}). So, we may apply a similar method and we conclude that $G^m$ is $0$-definable of finite width. Since $G/G^n$ is pseudofinite of finite exponent it is locally finite by Lemma  \ref{res-burn} and since it is finitely generated it must be finite. \qed

\bigskip
One may think to apply Proposition \ref{prop-def} to deduce immediately Proposition \ref{prop-fg-derived}. However, the problem is related to the following question. 

\frag Is a  $d$-generated pseudofinite group  pseudo-($d$-generated finite groups)? 
\efrag

\par We will use the following notation throughout the rest of this  section. Let $G$ be an infinite finitely generated pseudofinite group. Assume that $G$ is generated by $g_{1},\cdots,g_{d}$. By Frayne's theorem, there is an ultrapower   $\prod_I^{\U} F_{i}$ into which $G$ elementarily embeds.  Using this elementary embedding, we identify $g_{k}$, $1\leq k\leq d$, with $[f_{ki}]_{\U}$ with $f_{ki}\in F_{j}$. So, $G$ is isomorphic to the subgroup $\<[f_{1i}]_{\U},\cdots,[f_{di}]_{\U}\>$ of $\prod_{I}^{\U} F_{j}$.
\medskip

\prop \label{pseudo-type-fini} Let $G$ be a finitely generated pseudofinite group and suppose that
$G$ satisfies one of the following conditions. 
\begin{enumerate}
\item  $G$ is of finite exponent, or
\item {\rm (Kh\'elif)} $G$ is soluble, or
\item $G$ is soluble-by-(finite exponent), or
\item $G$  is pseudo-(finite linear of degree $n$ in characteristic zero), or
\item $G$ is simple, or 
\item $G$ is hyperbolic. 
\end{enumerate} 
Then such a group $G$ is finite. 
\eprop
\proof $(1)$  $G$ is locally finite by Lemma \ref{res-burn} and thus finite as it is finitely generated.  

\par $(2)$ Since $G$ is soluble, $G^{(n)}=1$ for some $n$ and thus $G$ is finite by Proposition \ref{prop-fg-derived}. 

 \par $(3)$ Assume  that $G$ is soluble-by-exponent $n$. Hence $G$ satisfies a nontrivial identity and thus by Lemma \ref{A_{n}} there exists $s$ such that for all $m\geq s$, the alternating groups $A_{m}$ are not involved in $F_{j}$. By Proposition \ref{prop-fg-derived} $G^n$ is $0$-definable and soluble. Since $G/G^n$ is a finitely generated pseudofinite group of exponent $n$,  by (1), it is finite, so $G^n$ is again a finitely generated soluble pseudofinite group and so it is finite by $(2)$. Thus $G$ is finite as required. 
 
 \par $(4)$ Let $G \preceq L=\prod_I^\U F_i$, where each $F_i$ is finite and linear of degree $n$ over $\Bbb C$. By a result of C. Jordan  \cite[Theorem 9.2]{weh}, there exists a function $d(n)$ depending only on $n$ such that each $F_i$ has an abelian subgroup of index at most $d(n)$.  Hence $L$ is abelian-by-finite and since $G$ is a subgroup of $L$, $G$ is also abelian-by-finite.  By $(3)$, $G$ is finite.  

\par $(5)$ In this case one uses Wilson's classification of the simple
pseudofinite groups \cite{W} and in particular the fact that they are all linear. Since $G$ is finitely generated and linear it is residually finite by a result of Mal'cev. Since $G$ is simple, it must be finite.   

\par $(6)$ If $G$ is not cyclic-by-finite then the commutator subgroup has an infinite width. Thus $G$ is cyclic-by-finite and thus finite by (3).  \qed

\frag Is a  pseudofinite linear group of degree $n$, pseudo-(finite and linear of degree $n$) ? 
\efrag

\frag Are there finitely generated infinite residually finite groups $G$ which are pseudofinite?
\efrag

\frag {\it (Sabbagh)} Are there finitely generated infinite groups $G$ which are pseudofinite?
\efrag

\section{Free subsemigroup, superamenability}

We study in this section the existence of free subsemigroups of rank two  in pseudofinite groups and its link with superamenability. Recall  that a group is \textit{superamenable} if  for any nonempty subset $A$ of $G$,  there exists a left-invariant finitely additive mesure  $\mu : P(G) \rightarrow [0, \infty]$ such that $\mu(A)=1$. It is known that a group containing a free subsemigroup of rank two is not superamenable \cite[Proposition 12.3]{para}. Superamenability is a strong form of amenability which  was introduced by Rosenblatt in \cite{Rosen} who also conjectured that a group is superamenable if and only if it is amenable and  does not contain a free subsemigroup of rank two. This question was settled negatively by R. Grigorshuck in \cite{Grig}. In this section, we show in particular, that for   $\aleph_0$-saturated pseudofinite groups, superamenability is equivalent to the absence of free subsemigroups of rank two.

\thm \label{thm-semi-superam} Let $G$ be an $\aleph_0$-saturated pseudofinite group. Then  either $G$ contains a free subsemigroup of rank $2$ or $G$ is nilpotent-by-(uniformly locally finite).  \ethm

\par Before proving Theorem 4.1, we will state two  corollaries. 
\defn\emph{\cite[Definition 12.7]{para} }

$\bullet$ Let $G$ be a group and $S$ a finite generating set of $G$. We let $\gamma_S(n)$ to be the cardinal of the ball of radius $n$ in $G$ (for the word distance with respect to $S\cup S^{-1}$), namely $\vert B^G_{S\cup S^{-1}}(n)\vert$ (see Notation 3.2). 

$\bullet$ A group $G$ is said to be \textit{exponentially bounded} if for any finite subset $S \subseteq G$, and any $b>1$, there is some $n_0 \in \Bbb N$ such that $\gamma_S(n) <b^n$ whenever $n>n_0$. 

\edefn

\kor \label{kor-semi-superam} Let $G$ be an $\aleph_0$-saturated pseudofinite group. Then the following properties are equivalent.

$(1)$ $G$ is superamenable. 

$(2)$  $G$ has no  free subsemigroup of rank $2$.  

$(3)$  $G$ is nilpotent-by-(uniformly locally finite).  

$(4)$ $G$ is nilpotent-by-(locally finite). 

$(5)$ Every finitely generated subgroup of $G$ is nilpotent-by-finite.  

$(6)$ $G$ is exponentially bounded. 

\ekor

\proof $(1) \Rightarrow (2)$.  This is exactly the statement of \cite[Proposition 12.3]{para}.

 $(2) \Rightarrow (3)$.  This is exactly the statement of Theorem \ref{thm-semi-superam}.

 $(3)   \Rightarrow (4) \Rightarrow (5)$. Clear. 
 
 $(5)   \Rightarrow (6) \Rightarrow (1)$. This  follows from  \cite{para}(see page 198 for more details).\qed

 \kor An infinite finitely generated pseudofinite group has approximation of degree $n$ to $M_2$ for every $n \in \Bbb N$. 
 
 \ekor
 
 \proof Suppose not and let $G \preceq L=\prod_I^\U G_i$, where each $G_i$ is finite. Then $L$ doesn't contain a free subsemigroup of rank $2$ and since it is $\aleph_0$-saturated, it is nilpotent-by-(uniformly locally finite). Therefore $G$ is nilpotent-by-finite and thus finite by Proposition \ref{pseudo-type-fini}; a contradiction. \qed

\bigskip
The rest of the section is devoted to the proof of  Theorem \ref{thm-semi-superam}. For $a,b \in G$, we let $H_{a,b}=\<a^{b^n}| n \in \mathbb Z\>$ and $H_{a,b}'$ its derived subgroup.  Let us recall the following definition (see \cite{P-milnor}, \cite{P}, \cite{mil}).
\defn  A nontrivial word $t(x,y)$ in $x, y$ is a {\it $N$-Milnor word} 
of {\it degree} $\leq \ell$ if it can be put in the form  $y x^{m_1} y^{-1} ... y^{\ell} x^{m_{\ell}} y^{-\ell} .u=1,$ where $u\in H_{x,y}'$, $\ell\geq 1$, $gcd(m_1,...,m_{\ell})=1$ (some of the $m_{i}$'s are allowed to take the value $0$)
and $\sum_{i=1}^{\ell} \vert m_{i}\vert\leq N$. 

\smallskip
\par A group $G$ is {\it locally N-Milnor} ({\it of degree $\leq \ell$}) if for all $a, b$ in $G$ there is a nontrivial $N$-Milnor word $t(x,y)$ (of degree $\leq \ell$) such that $t(a,b)=1$.

 \edefn
 It is straightforward that a group $G$ which contains the free group $F_{2}$, cannot be locally $N$-Milnor. 
\par Any nilpotent-by-finite group is locally $1$-Milnor. More generally one has the following property.
\lmm \emph{(\cite{Rosen} Lemma 4.8.)}

\label{lem-hab} Let $G$ be a group without free subsemigroup of rank $2$. Then for any  $a,b \in G$, the subgroup  $H_{a,b}$  is finitely generated, and $G$ is locally $1$-Milnor.
 \qed
\elmm
 \par A finitely generated linear group which is locally $N$-Milnor  is nilpotent-by-finite (see \cite{P} Corollary 2.3).

\medskip
\par
\noindent {\bf Example:} Let $p$ be a prime number and $C_{p}$ (respectively $C_{p^n}$) be the cyclic group of order $p$ (respectively $p^n$).  Then the finite metabelian groups $C_{p} wr C_{p^n}$, $n\in \omega-\{0\}$
do not satisfy an identity of the form $t(x,y)=1$, where $t(x,y)$ is a Milnor word of degree $<p^n$ (see \cite{P-milnor} Lemma 7).

\smallskip
\par On Milnor words, we will use the following theorem stated to G. Traustason (\cite{mil}). The key fact on these words is that the varieties of groups they define have 
the property that any finitely generated metabelian group in the variety is nilpotent-by-finite (\cite{BM} Theorem A). 
\medskip
\par To a Milnor word $t(x,y):=y x^{m_1} y^{-1} ... y^{\ell} x^{m_{\ell}} y^{-\ell} .u$, $u\in H_{x,y}'$, one associates a polynomial $q_{t}$ of $\Z[X]$ as follows: $q_{t}[X]=\sum_{i=1}^{\ell} m_{i}.X^i$ (see \cite{P-milnor}, \cite{mil}).
\thm {\rm(See Theorem 3.19 in \cite{mil})} Given a finite number of Milnor words $t_{i},\;i\in I$ and their associated polynomials $q_{t_{i}}, \,i\in I$,  there exist positive integers $c(q)$ and $e(q)$ only depending on $q:=\prod_{i\in I} q_{t_{i}}$, such that a finite group $G$ satisfying $\bigvee_{i\in I} t_{i}=1$, is nilpotent of class $\leq c(q)$-by-exponent dividing $e(q)$.\qed
\ethm
\par  Note that we can express by a universal sentence the property that a group $G$ is nilpotent of class $\leq c(q)$-by-exponent dividing $e(q)$.
So we can deduce the following. 
\kor Let $G$ be a group approximable by a class of finite groups which are locally   $N$-Milnor  of degree $\leq \ell$. Then $G$ is nilpotent-by-(uniformly locally finite).
\ekor
\pr By  Proposition \ref{prop-equiva}, $G\leq L=\prod_{i\in I}^{\U} F_{i}$, where each $F_{i}$ is a finite locally   $N$-Milnor group of degree $\leq \ell$. Hence for any $i$,   there is a finite disjunction $\bigvee_{j\in J_i} t_{j}(x,y)=1$, $J_i$ finite, where $t_{j}$ is a $N$-Milnor word of degree $\leq \ell$, such that  $F_{i}$ satisfies $\bigvee_{j\in J_i} t_{j}(x,y)=1$. Let $q_i:=\prod_{j\in J} q_{t_{j}}$. By the theorem above, there exist positive integers $c(q_i)$ and $e(q_i)$, such that $F_{i}$ is nilpotent of class $\leq c(q_i)$-by-exponent dividing $e(q_i)$. 
Since the degree of each $q_{t_{j}}$ is bounded by $\ell$ and their coefficients are bounded in absolute value by $N$, there are a finite number of such polynomials. Let $Q$ denote the set of all possible products of such polynomials. Let $c_{max}:=max\{c(q): q\in Q\}$ and $e_{max}:=\prod_{q\in Q} e(q)$. So for each $i\in I$, we have that $F_{i}^{e_{max}}$ is nilpotent of class $\leq c_{max}$.  Set $N=\prod_{i\in I}^{\U} F_{i}^{e_{max}}$. Then $N$ is nilpotent and since $L^{e_{max}} \leq \prod_{i\in I}^{\U} F_{i}^{e_{max}}$ we conclude that $L/N$  is of finite exponent. Since $L/N$ is pseudofinite, by Lemma \ref{res-burn} $L/N$ is uniformly locally finite. Thus $L$ is nilpotent-by-locally finite as well as $G$.   \qed

\bigskip
\par
 \noindent {\bf Proof of Theorem \ref{thm-semi-superam}}
\pr   Let $G$ be an $\aleph_{0}$-saturated pseudofinite group not containing the free subsemigroup of rank $2$. Then, by Lemma \ref{lem-dicho}, it satisfies a finite disjonction of positive identities. In particular there exists $\ell$ such that  
it is approximable by a class of finite groups locally $1$-Milnor of degree $\leq \ell$ and so we may apply the preceding corollary.\qed
\bigskip
 \kor An $\aleph_{0}$-saturated locally $N$-Milnor pseudofinite group
 is nilpotent-by-(uniformly locally finite).
 \ekor
 \pr It is proven in the same way as the above theorem, using a similar argument as in Lemma \ref{lem-dicho} to show that such group satisfies a finite disjonction of identities of the form $t_{i}(x,y)=1$, where $t_{i}(x,y)$ is a $N$-Milnor word. And so again we can find a bound on the degrees of the corresponding Milnor words. \qed
\bigskip
\par

\noindent {\bf Example:} Y. de Cornulier and A. Mann have shown that if one takes the non Milnor word $[[x,y],[z,t]]^q$, then there is a residually finite $2$-generated group satisfying the identity $[[x,y],[z,t]]^q=1$, which is not soluble-by-finite (\cite{CM}). They exhibit a family of finite soluble groups $R_{n}$ generated by two elements, of solubility length $n$  and satisfying the identity $[[x,y],[z,t]]^q=1$. 
\par Let us recall their construction. On one hand they use an embedding theorem due to B.H. Neumann and H. Neumann in wreath products (\cite{NN}) and on the other hand a result of Razmyslov (\cite{V} chapter 4) that for each prime power $q\geq 4$ there exist a finite group $B_{r}$ generated by $r$ elements, of exponent $q$ and solubility length $n:=\lfloor log_{2}(r)\rfloor$. By \cite{NN}, $B_{r}$ embeds in a two generated subgroup $R_{n}$ of $(B_{r}Wr C_{p^k}) Wr C_{p^k}$, for some $k$ sufficiently large. (The number $k$ is chosen 
such that $p^k \geq 4r-1$.)
So $R_{n}$ is a $2$-generated $p$-group satisfying the identity $[[x,y],[z,t]]^q=1$.
\par In particular, we have an example of an $\aleph_{0}$-saturated pseudofinite group $L$ not containing $F_{2}$ and not soluble-by-finite (with $\phi_{R}(L)=L$) (see Proposition \ref{segal-res}).
 Take $L=\prod_{\Bbb N}^{\U} R_{n}$. 

\bigskip

\section{Free subgroups, amenability}

As we have seen  in the previous section, the absence of free subsemigroups of rank $2$ in pseudofinite ($\aleph_0$-saturated) groups implies superamenability. In this section,  we are interested in the similar problem with  free subgroups of rank $2$. However, as the next proposition shows, the problem is connected to some strong properties that residually finite  groups must  satisfy.  We will be interested in this section, more particularly,  in the problem of  amenability of pseudofinite groups. 
Then in the  next section, we shall give some alternatives under stronger hypotheses. 
\par Recall that a group is said to be \textit{amenable} if there exists a finitely additive left-invariant  measure $\mu : \P(G) \rightarrow [0, 1]$ such that  $\mu(G)=1$. Note that
there are many definitions of amenable groups in the literature (see for instance Theorem 10.11 \cite{para}).  
\par M. Bozejko \cite{Boz}  and G. Keller \cite{Kel} called a group $G$ \textit{uniformly amenable} if there exists a function $\alpha : [0,1] \times \N \rightarrow \N$ such that for any finite subset $A$ of $G$ and every $\epsilon \in [0,1]$ there is a finite subset $V$ of $G$ such that $|V|\leq \alpha (\epsilon, |A|)$ and $|AV|<(1+\epsilon)|V|$. 
By using the equivalent definition of amenability with F\o lner sequences, we have that an uniformly amenable group is amenable.

\thm \label{prop-equiv} The following properties are equivalent.  

\begin{enumerate}

\item  Every $\aleph_0$-saturated pseudofinite group either contains a free nonabelian group or it is amenable. 

\item Every ultraproduct of finite groups either contains a free nonabelian group or it is amenable. 

\item  Every finitely generated residually finite group satisfying a nontrivial identity is amenable. 

\item Every finitely generated residually finite group satisfying a nontrivial identity is uniformly amenable.

\end{enumerate}
\ethm
\par G. Keller showed that a group $G$ is uniformly amenable if and only if all its ultrapowers are amenable. Later J. Wysoczanski \cite{wy} gives a more simple combinatorial proof. However, the notion which is behind this is the saturation property. 

\bem Let $\sigma_{p,n,f}$ be the following sentence with $(p,n)\in \N^2$ and $f:\N^2\rightarrow \N$:
$$\forall a_{1}\cdots\forall a_{n}\exists y_{1}\cdots\exists y_{f(p,n)}\;
p.\vert\{a_{i}.y_{j}: 1\leq i\leq n; 1\leq j\leq f(p,n)\}\vert< (p+1).f(p,n).$$
Then we see that $G$ is uniformly amenable iff there exists a function $f:\N^2\rightarrow \N$ such that $G\models \sigma_{p,n,f}$ for any $(p,n)\in \N^2$. In particular being uniformly amenable is elementary, that is it is a property preserved by elementary equivalence.  \ebem

\prop \label{prop-aleph-ame}An $\aleph_0$-saturated group is amenable if and only if it is uniformly amenable. \eprop

\proof Suppose that $G$ is $\aleph_0$-saturated and  amenable. Then for any finite subset $A$ of $G$ and every $\epsilon \in [0,1]$ there is a finite subset $V$ of $G$ such that  $|AV|<(1+\epsilon)|V|$.

Let $A=\{a_1, \dots, a_n\}$ and $\epsilon \in[0,1]$. We may assume without loss of generality that $\epsilon=1/p$ for some $p \in \N$. Then 
$$
G \models \bigvee_{m \in \N} \exists x_1 \dots \exists x_m (p|A.\{x_1, \dots, x_m\}|<(p+1).m),
$$
and by $\aleph_{0}$-saturation
$$
G \models \bigvee_{1 \leq m \leq r} \exists x_1 \dots \exists x_m (p|A.\{x_1, \dots, x_m\}|<(p+1).m.
$$
By setting $\alpha(\epsilon, n)=r$, we get the uniform bound. Thus $G$ is uniformly amenable. \qed
\kor {\rm (\cite{Kel}, \cite{wy})} A group is uniformly amenable if and only if all its nonprincipal ultrapowers are amenable iff one of its nonprincipal ultrapowers is amenable. \qed
\ekor

\par In the proof of Theorem \ref{prop-equiv}, we will use the fact that the class of amenable groups is closed under various operations (see \cite{para} Theorem 10.4) and in particular a group is amenable iff its finitely generated subgroups are. Recall that no amenable group contains the free subgroup of rank $2$ (\cite{para} Corollary 1.11).
\par We will need the following simple lemma (\cite{Kel} Theorem 4.5). 

\lmm  \label{unif-sub} A subgroup of an uniformly amenable group is uniformly amenable.  \qed
\elmm
\bigskip
\noindent
\textbf{Proof of Theorem \ref{prop-equiv}.} 

$(1) \Rightarrow (2)$. An ultraproduct of finite groups  is $\aleph_0$-saturated and pseudofinite, the conclusion follows. 

$(2) \Rightarrow (3)$.  Let $G$ be a finitely generated residually finite group satisfying a nontrivial identity $t=1$.  Then $G$ embeds into an ultraproduct $K$ of finite groups which satisfies a nontrivial identity (see Proposition 2.2 and note that $G$ is residually $\C$ with $\C$ the class of finite groups satisfying $t=1$). By $(2)$,  $K$ is amenable and thus $G$ is amenable.

$(3) \Rightarrow (1)$. Let $G$ be an $\aleph_0$-saturated pseudofinite group and suppose that $G$ has no free nonabelian subgroup.  Let $K$ be an ultraproduct of finite groups  such that $G \preceq K$. Since $G$ is $\aleph_0$-saturated, $G$ satisfies a nontrivial identity by Lemma \ref{lem-dicho},  as well as $K$.  It is sufficient to show that every finitely generated subgroup of $K$ is amenable.   Let $L$ be a finitely generated subgroup of $K$.  Let $\C$ be the class of finite groups  satisfying  the identity satisfied by $K$. Then $L$ is approximable by $\C$, and since $\C$ is a pseudovariety, by Proposition \ref{prop-equiva} $L$ is a direct limit of fully residually-$\C$ groups. Hence $L$ is a direct limit of residually finite groups satisfying a nontrivial identity. By our hypothesis such groups are amenable as well as their direct limit $L$.

Clearly $(4) \Rightarrow (3)$ and it remains to show $(3) \Rightarrow (4)$.  Let $L$ be a finitely generated residually finite group satisfying a nontrivial identity, so $L$ is residually $\C$, where $\C$ is a class of finite groups satisfying a nontrivial identity. By Proposition 2.2, $L$ is approximable by $\C$, namely
embeds in an ultraproduct $K$ of elements of $\C$. By $(1)$, $K$ is amenable. Since $K$ is $\aleph_0$-saturated, it is uniformly amenable by Proposition \ref{prop-aleph-ame} as well as $L$ by Lemma \ref{unif-sub}. 

\qed

\bigskip

\par As recalled above, a group containing a nonabelian free group cannot be amenable. Von Neumann  and Day asked for the converse, namely whether every non-amenable group contains a nonabelian free group. This was answered negatively by Ol'shanskii \cite{Ol}, Adyan \cite{Ad} and Gromov \cite{Gro}.  However a positive answer can be provided for some classes of groups as the class of linear groups. One may ask if the question has a positive answer in the class of residually finite groups. This was answered negatively by Ershov \cite{Er}. Other examples  of non-amenable residually finite groups without nonabelian free subgroups were constructed by Osin in \cite{Os}. 

\frag \cite[Question 14]{CM} Does there exist a non-amenable finitely generated residually finite group satisfying a nontrivial identity? 
\efrag

\defn Let $\mathcal G=(G_i)_{i \in I}$ be a family of groups and $\U$ an ultrafilter over $I$. We say that $\mathcal G$ is \textit{uniformly amenable} relative to $\U$ if the following condition holds.   There exists a function $\alpha : [0,1] \times \N \rightarrow \N$ such that for any $n \in \N$  and every $\epsilon \in [0,1]$, there exists $U \in \U$ such that for any $i \in U$, for any finite subset $A $ in $G_i$ with $|A|=n$,  there is a finite subset $V$ of $G_i$ such that $|V|\leq \alpha (\epsilon, |A|)$ and $|AV|<(1+\epsilon)|V|$.
\edefn

A proof similar to that of Proposition \ref{prop-aleph-ame} yields. 

\prop Let $\mathcal G=(G_i)_{i \in I}$ be a family of groups and $\U$ an ultrafilter over $I$. Then $\prod_I^\U G_i$ is amenable if and only if $\mathcal G$ is uniformly amenable relative to $\U$. \qed
\eprop

\frag Is a pseudofinite amenable group uniformly amenable? 
\efrag

In \cite{HPP}, the notion of \textit{definably amenable} groups was introduced. A group is said to be \textit{definably amenable}  if there exists a finitely additive left-invariant measure $\mu : \D(G) \rightarrow [0,1]$  with $\mu(G)=1$; where $\D(G)$ is the boolean algebra of definable subsets of $G$.   They pointed out that there are definably amenable groups which are not amenable as $SO_{3}(\R)$ and also groups that are not definable amenable as $SL_{2}(\R)$.

\par The following proposition gives natural examples of definably amenable groups (and again shows that there are definably amenable but non amenable groups).

\prop \label{pseudo-def-amenable}A pseudofinite group is definably amenable. 
\eprop

We will show in fact the next more general proposition. Let us first give a definition which is borrowed from non standard analysis.

\defn Let $I$ be a set and $(G_i)_{i \in I}$ a family of groups, $\U$ and ultrafilter on $I$. A subset $A \subseteq \prod_I^\U G_i$ is said to be \textit{internal} if there exists $(A_i)_{i \in I}$, $A_i \subseteq G_i$,  such that $A= \prod_I^\U A_i$. 
\edefn

We see that every definable subset is internal and that the set of internal subsets forms a left-invariant  boolean algebra.  Recall that a measure on a boolean algebra $\B$,  is said \textit{$\sigma$-additive} if $\mu(\bigcup_{n \in \Bbb N}A_i)=\sum_{i \in \Bbb N}\mu(A_i)$ whenever $A_i \cap A_j =\emptyset$ and $\bigcup_{i \in \Bbb N}A_i \in \B$.  Given a boolean algebra $\B$, we denote by $\bar \B$ the $\sigma$-algebra generated by $\B$. 

\prop \label{prop-amenable-ultra} Let $(G_i)_{i \in I}$ be a family of amenable groups, $\U$ an ultrafilter on $I$.  Let $\B$ be the boolean algebra of internal subsets of $L=\prod_I^\U G_i $. Then there exists a finitely additive measure  $\mu : \P(L) \rightarrow [0,1] $, $\mu(L)=1$,  whose restriction to $\bar \B$ is $\sigma$-additive and left-invariant.  
\eprop

\pr   If  $A \in \mathcal B$ and  $A= \prod_I^{\U}A_i$,  we define the measure of  $A$ by 
$$
\mu(A)= \lim_{\mathcal U}\mu_i(A_i),
$$
 where each $\mu_i$  is a left-invariant probability measure on $G_i$.  It is not difficult to see that $\mu$ is a left-invariant  finitely additive measure defined on $\B$.

 Let us show that $\mu$ is $\sigma$-additive.  It is sufficient to show that if  $(B_i| i \in \mathbb N)$  is a sequence of  $\B$, such that  $B_{i+1} \subseteq B_i$ and  $\cap B_i =\emptyset$, then $\lim_{i \rightarrow \infty}\mu(B_i)=0$.

But by the saturation of the ultraproduct,  we get  $B_1 \cap \dots \cap B_n =\emptyset$.  Hence  $\mu(B_i)=0$  for  $i \geq n$  and so  $\lim_{i \rightarrow \infty}\mu(B_i)=0$.

By  Carath\'eodory's theorem, $\mu$ can be extended  to a  $\sigma$-additive $\bar \mu$  defined over  $\bar{\mathcal B}$.  It is not difficult to see that $\bar \mu$ is still left-invariant on $\bar \B$.   By a theorem of  Horn and  Tarski \cite{para}, $\bar \mu$  can be extended  to a finitely additive measure defined on  $\P(L)$. 
\qed

\bigskip
\noindent
\textbf{Proof of Proposition  \ref{pseudo-def-amenable}.} Let $G$ be a pseudofinite group. Then $G \preceq \prod_I ^\U G_i$, where each $G_i$ is a finite group. By Proposition \ref{prop-amenable-ultra}, there exists a left-invariant probability measure $\mu$ defined on definable subsets of $L$. For every definable subset $X$ of $G$,  definable by a formula $\phi(x)$, we take $\mu (X)=\mu(\phi(L))$. Then this defines a left-invariant probability measure $\mu$  on definable subsets of $G$. \qed

\section{Free subgroups, alternatives}

We study in this section existence of free subgroups in pseudofinite groups under strong hypotheses. Recall that the \textit{Pr\"ufer} rank of a  group $G$ is the least integer $n$ such that every finitely generated subgroup of $G$ can be generated by $n$ elements.  S. Black \cite{B} has considered famillies $\C$ of finite groups of bounded Pr\"ufer rank and showed that a finitely generated residually-$\C$ group $G$ either contains a free nonabelian group or it is nilpotent-by-abelian-by-finite. 

S. Black (Theorem A in \cite{B2})  also showed the following \textit{finitary Tits' alternative}:  there exists a function $d(n,r)$ such that if $G$ is a finite group of \textit{Pr\"ufer} rank $r$ then either $G$ contains an approximation of degree $n$ to $F_2$ or $G$ has a soluble subgroup whose derived length and index in $G$ are at most $d(n,r)$. Moreover, in this case, there exists $c=c(n,r)$, $\ell=\ell(n,r)$ such that $G$ is nilpotent of class at most $c$-by-abelian-by-index-at-most-$\ell$. 

First, we reformulate the above result in the context of pseudofinite groups. 

\thm \label{thm-alter-bounded-rank}Let $G$ be an $\aleph_0$-saturated pseudo-(finite of bounded Pr\"ufer rank) group. Then either $G$ contains a nonabelian free group or $G$ is nilpotent-by-abelian-by-finite. 
\ethm 
\kor An $\aleph_0$-saturated pseudo-(finite of bounded Pr\"ufer rank) group either contains a nonabelian free group or is uniformly amenable.\qed
\ekor

\bem The above theorem is equivalent to the previous mentioned result of Black.  One direction is clear: the finitary Tits' alternative implies that if $L=\prod_I^{\U} F_i$, where $L$ satisfies a nontrivial identity then for almost $i$, $G_i$ is nilpotent-by-abelian-by-finite with an uniform bound on the nilpotency classes  and the indices; hence $L$ is nilpotent-by-abelian-by-finite. For the converse  let $(G_n|n \in \mathbb N)$ be the sequence of all finite groups of Pr\"ufer rank $r$ and without an approximation of degree $n$ to $F_2$.  Suppose that for any $d \in \mathbb N$ there exists $G_{n_d}$ such that $G_{n_d}$ is not nilpotent of class at most $d$-by-abelian-by-index-at-most-$d$.  Hence we get  an infinite sequence $H_d=G_{n_d}$. Let $\U$ be any nonprincipal ultrafilter on $\mathbb N$ and set $L=\prod_{\mathbb N}^{\U}H_d$. Then $L$ satisfies a nontrivial identity  and thus $L$ is nilpotent of class at most $p$-by-abelian-by-index-at-most-$p$ for some $p \in \mathbb N$.  Since this last property can be expressed by a first order sentence,  there exists $U \in \U$ such that for any $d \in U$,  $H_d$ is nilpotent of class at most $p$-by-abelian-by-index-at-most-$p$.  Hence for $d\geq p$ and $d \in U$,  $H_d$ is nilpotent of class at most $d$-by-abelian-by-index-at-most-$d$; which is  clearly a contradiction.   
\ebem

We give here the proof of Theorem \ref{thm-alter-bounded-rank} from the pseudofinite groups viewpoint. As in \cite{B2}, one reduces first the problem to finite soluble groups, using a result of A. Shalev (\cite{Sh}) and then one uses a result of D. Segal on residually (finite soluble groups) (\cite{Se-rank}).  
\par We first note that the following alternative holds for simple pseudofinite groups.
\lmm \label{lem-jones} Let $G$ be a nonprincipal ultraproduct of finite simple groups. Then either $G$ contains a nonabelian free group, or $G$ is finite.
\elmm
\pr Assume that $G$ doesn't contain a nonabelian free group. Then, since $G$ is $\aleph_{0}$-saturated, by Proposition \ref{lem-dicho}, $G$ satisfies a nontrivial identity. By a theorem of Jones that a proper variety of groups only contains finitely many finite nonabelian simple groups (\cite{J}), we can bound the cardinality of the finite simple groups appearing in that ultraproduct, which contradicts the fact that $G$ is infinite.
\qed
\kor  \label{jones-simple}  Let $G$ be a simple pseudofinite group. Then either $G$ contains a nonabelian free group, or $G$ is finite.
\ekor
\pr By a result of J. Wilson (\cite{W0}) and its strenghtening (\cite{EJMR}), $G$ is isomorphic to a non principal ultraproduct of finite simple groups (of fixed Lie type). Then we apply the above lemma.
\qed
\medskip
\par Recall that a group is said \textit{quasi-linear} if it is embeddable in a finite direct product of linear groups.  We say that a function is  \textit{$r$-bounded} if it is bounded in terms of $r$ only. 

\prop Let $G$ be a semi-simple pseudo-(finite of bounded Pr\"ufer rank) group. Then $G$ has a quasi-linear subgroup of finite index.
\eprop

\proof  Let $G\prec L$, where $L= \prod_I^\U G_i$, $\U$ a nonprincipal ultrafilter on $I$ and each $G_i$ is finite of bounded \textit{Pr\"ufer} rank, $i\in I$.  Since $\phi_{R}(G)=\{1\}$, we have that $\phi_{R}(L)=\{1\}$ and so on an element $U$ of $\U$, each $G_{i}$ is semi-simple. Using Proposition 3.6 in \cite{Sh}, for $i\in U$, there exists a characteristic subgroup $G_{1i}$ of $G_{i}$ such that  $|G_i/G_{1i}|$ is $r$-bounded, say of cardinality $\leq f(r)$ and  $G_{1i} \cong S_{1i} \times \dots \times S_{ki} $ where $1 \leq k \leq g(r)$ and each $S_{1i}$ is a simple pseudo-finite group of Lie type of $r$-bounded Lie rank $n_j$ over the finite field $F_{p_j^{e_j}}$ where $e_j$ is $r$-bounded and $1\leq j\leq r$, for each $i\in U$. 
\par We have $1 \,\triangleleft \prod_{I}^{\U} G_{1i}\; \triangleleft\; \prod_{I}^\U G_i=L$ and the subgroup $L_{0}:=\prod_{I}^{\U} G_{1i}$ is of finite index in $L$ since $\vert \prod_{I}^{\U} G_{i}/\prod_{I}^{\U} G_{1i}\vert \leq f(r)$. Moreover,  $\prod_{I}^{\U} G_{1i}\cong \prod_I^\U  (S_{1i} \times \dots \times S_{ki}) \cong  (\prod_I^\U S_{1i}) \times \dots \times (\prod_I^\U S_{ki})$ and each factor is a simple linear group.  
Since $G$ embeds in $L$ and $G\cap L_{0}$ is a subgroup of finite index in $G$ which embeds in a quasi-linear group, $G$ has a quasi-linear group of finite index. \qed

\kor \label{semi-simple} Let $G$ be a pseudo-(finite of bounded Pr\"ufer rank) group. Then $G/\phi_R(G)$ has a quasi-linear subgroup of finite index.

\ekor
\pr This follows from the previous  proposition and Lemma \ref{radical}.
\qed

\prop \label{segal-res} Let $G$ be a pseudo-(finite of bounded Pr\"ufer rank) group satisfying a nontrivial identity and such that $G=\phi_R(G)$. Then $G$ is nilpotent-by-abelian-by-finite. 
\eprop

\proof Let $G \preceq L= \prod_I^\U G_i$, where each $G_i$ is finite of  Pr\"ufer rank $\leq r$. 
Since $G=\phi_{R}(G)$ and $G\preceq L$, without loss of generality we may assume that each $G_i$ is soluble. Similarly since $G$ satisfies a nontrivial identity, $L$ satisfies a nontrivial identity, say $t=1$. Hence, 
 w.l.o.g. each $G_i$ satisfies the same nontrivial identity $t=1$. 
\par Let $F=\<x_1, \dots, x_r|\>$ be the free group on $\{x_1, \dots, x_r\}$ and for each $i \in U$ let $S_i=\{s_{1i}, \dots, s_{ri}\}$ be a finite generating set of $G_i$. Let $\varphi_i : F \rightarrow G_i$ be the natural homomorphism which sends $x_j$ to $s_{ji}$. Let $H=F/ \bigcap_{i \in U}\ker(\varphi_i)$.  Then $H$ is  a residually (finite soluble of bounded rank) group, satisfying a nontrivial identity $t=1$. 
By a result of Segal (see Theorem page 2 in \cite{Se-rank}), $H$ has a nilpotent normal subgroup $N$ such that $H/N$ is quasi-linear.  Since $H/N$ doesn't contain $F_{2}$, by Tits alternative for linear groups, $H/N$ is soluble-by-finite. Hence $H$ has a soluble normal subgroup $K$ of finite index. Again by the same theorem of Segal (\cite{Se-rank}), $K$ is nilpotent-by-abelian-by finite, so is $H$. Hence there exists $c$ and $f$ such that each $F/\ker(\varphi_i)\cong G_i$ is (nilpotent of class at most $c$)-by-abelian-by-finite index $f$. So, $[L^f,L^f]$ is nilpotent of class at most $c$.
\par Since $L$ does not contain $F_{2}$, by Proposition \ref{prop-def}, $L^f$ is of finite index in $L$, $L^f$ is $0$-definable and so we can express in a first-order way that $[L^f,L^f]$ is nilpotent of class at most $c$. These (first-order) properties transfer in $G$.
\qed

\bigskip
\noindent \textbf{Proof of Theorem \ref{thm-alter-bounded-rank}.} Let $G \preceq L= \prod_I^\U G_i$, where each $G_i$ is finite of rank $\leq r$. Suppose that $G$ contains no free nonabelian subgroup. Then $G$ satisfies a nontrivial identity by Lemma \ref{lem-dicho}, as well as $L$ and $L/\phi_R(L)$. 
By Corollary \ref{semi-simple} and Lemma \ref{lem-jones}, $G/\phi_{R}(G)$ is finite, say of cardinality $\leq f(r)$.
Since $G$ does not contain $F_{2}$, by Proposition \ref{prop-def}, $G^{f(r)}$ is $0$-definable and of finite index in $G$.
Applying Proposition \ref{segal-res} to $G^{f(r)}$, we get that $G^{f(r)}$ is nilpotent-by-abelian-by-finite. So the conclusion also applies to $G$.\qed


\medskip

We place ourselves now in a slightly more general context than Theorem \ref{thm-alter-bounded-rank}. 
\defn Let us say that a class $\C$ of finite groups is \textit{weakly of $r$-bounded rank} if for each element $G\in \C$, the index of the sockel of $G/rad(G)$ is $r$-bounded and $rad(G)$ has $r$-bounded rank.
\edefn

\par By the above result of A. Shalev (\cite{Sh}), a class of finite groups of $r$-bounded Pr\"ufer rank is weakly of $r$-bounded rank. 
\defn {\rm (See \cite{res}.)} A group $G$ has {\it finite $c$-dimension} if there is a bound on the chains of centralizers.
We will say that a class $\C$ of finite groups has {\it bounded $c$-dimension} if there is $d\in \N$ such that for each element $G\in \C$, the $c$-dimensions of $rad(G)$ and of the sockel of $G/rad(G)$ are $d$-bounded. (Note that a class of finite groups of bounded Pr\"ufer rank is of bounded $c$-dimension.)
\edefn
\lmm \label{bounded} Let $\C$ be a class of finite groups satisfying a nontrivial  identity. Suppose that for any $G\in \C$, $Soc(G/rad(G))$   is of $r$-bounded rank, or of $r$-bounded index in $G/rad(G)$ or of $r$-bounded $c$-dimension. Then $G/rad(G)$ is of bounded exponent depending only on $r$ and the identity. 
\elmm
\pr Recall that the sockel $Soc(G)$ of a group $G$ is the subgroup generated by  all   minimal normal nontrivial subgroups of $G$. In case $G$ is a finite group, then 
$Soc(G)$ is a direct sum of simple groups and is completely reducible (see \cite{S} 7.4.12).  
\par Let $S:=Soc(G/rad(G))$. Since a nontrivial identity can only be satisfied by finitely many finite simple groups (see \cite{J}), by hypothesis on the class $\C$, we have a bound on the cardinality of the simple groups appearing in $Soc(G/rad(G))$, for $G\in \C$. So if the index of $Soc(G/rad(G))$ in $G/rad(G)$ is $r$-bounded, then the exponent of $G/rad(G)$ is bounded in terms of $r$ and the identity only. 

In the other cases, we note the following. The centralizer of $S$ is trivial and so in order to show that some power of an element of $G/rad(G)$ is equal to $1$, it suffices to show that the corresponding inner automorphism on $S$ is the identity. 
\par Let $\bar g\in G/rad(G)$ and let $\alpha_{\bar g}$ be the conjugation by $\bar g$ in $G/rad(G)$. It induces a permutation of the copies of a given finite simple group appearing in $S$. So if the subgroups of $S$ generated by $\alpha_{{\bar g}^z}(\bar h)$, $z\in \Z$, are $r$-generated, or if the $c$-dimension of $S$ is $r$-bounded, we get the result.
\qed
\kor Let $G$ be a pseudo-(finite weakly of $r$-bounded rank) group satisfying a nontrivial identity. Then $G/\phi_R(G)$ is 
uniformly locally finite. 
\ekor

\proof Let $G \preceq L= \prod_I^\U G_i$, where each $G_i$ is finite. 
Let $\C:=\{G_{i}\,:\;i\in I\}$, then it satisfies the hypothesis of the lemma above. So, $L/\phi_{R}(L)$ is of bounded exponent. It transfers to $G/\phi_{R}(G)$. We conclude by applying Lemma \ref{res-burn}.\qed

\thm \label{thm-alter-weak} Let $G$ be an $\aleph_0$-saturated pseudo-(finite weakly of bounded rank) group.  Then either $G$ contains a nonabelian free group or $G$ is nilpotent-by-abelian-by-(uniformly locally finite). 
\ethm

\pr By applying Proposition \ref{segal-res} to $\phi_{R}(G)$, we get that $\phi_{R}(G)$ is nilpotent-by-abelian-by-finite. By the above Corollary, $G/\phi_{R}(G)$ is uniformly locally finite. So, $G$ is nilpotent-by-abelian-by-uniformly locally finite.\qed
\medskip

\prop Let $\C$ be a class of finite groups of bounded $c$-dimension and suppose $G$ is a pseudo-$\C$ group satisfying a nontrivial identity. Then $G$ is soluble-by-(uniformly locally finite). 
\eprop
\proof Let $G \preceq L= \prod_I^\U G_i$, where each $G_i$ is finite. 
Then by hypothesis there is $d\in \N$ such that the $c$-dimension of each $\phi_{R}(G_{i})$ is $d$-bounded (as well as the $c$-dimensions of $Soc(G_{i}/rad(G_{i}))$) and so by a result of E. Khukhro (see Theorem 2 in \cite{res}), the derived length of $\phi_{R}(G_{i})$ is $d$-bounded. So, $\phi_{R}(L)$ is soluble as well as $\phi_{R}(G)$.
By Lemma \ref{bounded}, the $G_i/\phi_R(G_i)$ are of an unfirom bounded exponent and so $L/\phi_{R}(L)$ is of finite exponent as well as $G/\phi_{R}(G)$. By Lemma \ref{res-burn}, $G/\phi_{R}(G)$ is uniformly locally finite.\qed

\kor \label{thm-finite-dim} Let $G$ be an $\aleph_0$-saturated pseudo-(finite of bounded $c$-dimension) group.  Then either $G$ contains a nonabelian free group or $G$ is soluble-by-(uniformly locally finite). 

\qed
\ekor

\end{document}